\documentclass[11pt,a4paper]{article}
\usepackage{amssymb,amsmath,amsfonts,amscd,amsthm,amsxtra,latexsym}
\usepackage[T1]{fontenc}    
\usepackage{stmaryrd}
\theoremstyle{plain}

\newtheorem{corollary}{Corollary}

\theoremstyle{remark}
\newtheorem{definition}{Definition}
\newtheorem{remark}{Remark}

\,
\newcommand{\SO}{\mathrm{SO}}

\newcommand{\SU}{\mathrm{SU}}

\newcommand{\Sp}{\mathrm{Sp}}
\newcommand{\so}{\mathfrak{so}}
\newcommand{\su}{\mathfrak{su}}

\renewcommand{\sp}{\mathfrak{sp}}

\newcommand{\frakc}{\mathfrak{c}}

\newcommand{\frakm}{\mathfrak{m}}
\newcommand{\frakn}{\mathfrak{n}}
\newcommand{\frakt}{\mathfrak{t}}
\newcommand{\Sym}{\mathrm{Sym}}

\newcommand{\g}[2]{\langle #1,#2\rangle}
\renewcommand{\gg}{\g{\,\cdot\,}{\,\cdot\,}}

\newcommand{\R}{\mathrm{I\!R}}

\newcommand{\N}{\mathbb{N}}

\newcommand{\frakg}{\mathfrak{g}}
\newcommand{\frakh}{\mathfrak{h}}

\newcommand{\dt}{\mathrm{dt}}

\newcommand{\diag}{\mathrm{diag}}

\renewcommand{\i}{\mathrm{i}}
\renewcommand{\j}{\mathrm{j}}
\renewcommand{\k}{\mathrm{k}}
\newcommand{\bbH}{\mathbb{H}}

\,

\newcommand{\rmS}{\mathrm{S}}
\newcommand{\rmU}{\mathrm{U}}

\newcommand{\scrR}{\mathcal{R}}

\newcommand{\jet}{\mathfrak{jet}}
\newcommand{\fraku}{\mathfrak{u}}

\newcommand{\Cl}{\mathrm{Cl}}

\title{An estimate for the Singer invariant via the Jet Isomorphism Theorem}
\author{Tillmann Jentsch}
\begin{document}\sloppy

\maketitle
\begin{abstract}
Recently examples of Riemannian homogeneous spaces with linear Jacobi relations were found. 
We calculate the Singer invariants of these spaces with the computer algebra program Maple and discuss the results by means of the Jet Isomorphism Theorem of Riemannian geometry.
\end{abstract}


\section{Overview}
A Riemannian space $M$ satisfies a linear Jacobi relation of order $k\geq 0$ if
the $k + 1$-th covariant derivative of the Jacobi operator along an arbitrary
unit-speed geodesic is a linear combination of  covariant derivatives of lower
order with coefficients independent of the geodesic. Examples are given by the g.o-spaces with {\em constant Jacobi osculating ranks} found in~\cite{A1,A2,A3,A4,G,GN,NT}.

For the following homogeneous spaces it was already shown that they have
constant Jacobi osculating ranks and, in particular, satisfy a linear Jacobi relation of order
$2$ and $4$\,, respectively: 

the six-dimensional complex flag manifold $M^6$ equipped with the normal
homogeneous metric, the exceptional normal homogeneous spaces of  positive sectional curvatures known as the Berger manifold $V_1$ and Wilkings example $V_3$ (with the standard metric)\,, and for the two-step nilpotent group $N^6$ known as Kaplans example of a g.o. space which is not naturally reductive.

According to the {\em Jet  Isomorphism Theorem} of Riemannian geometry, the
Singer invariant of a Riemannian homogeneous space $M$ can not exceed the constant $k$\,. Using the computer algebra program Maple, we will calculate the Singer invariants for the above examples explicitly and see that their Singer invariants are in fact much smaller.

\section{Linear  Jacobi relations}
Let $(M,\gg)$ be a pseudo-Riemannian manifold and $p\in M$\,.  We denote by $\nabla$
and $R$ the Levi-Civita connection and the curvature tensor. Recall that the Jacobi operator along some geodesic $\gamma$ is defined by $\scrR_\gamma(t;x,y) := R(x,\dot\gamma(t),\dot\gamma(t),y)$ for all $t\in \R$ and $x,y\in T_{\gamma(t)}M$\,. Further, let $\Sym^kV$ denotes the  $k$-th symmetric power of an arbitrary vector space $V$\,.
Thus the Jacobi operator is a section of $\Sym^2 TM^*$ along $\gamma$\,. 

Furthermore, let $\nabla^{k)}R$ denote the $k$-th iterated covariant derivative of the curvature tensor. Then the $k$-th covariant derivative of $\scrR_{\gamma}$ is given by
\begin{equation}\label{eq:k-ter_Jacobi_on_a_geodesic}
\scrR^{k)}_\gamma(x,y) := \frac{\nabla^{k}}{\dt^k} \scrR_\gamma(x,y) = \nabla^{k)}_{\dot \gamma,\cdots,\dot \gamma}R(\dot \gamma,x ,y ,\dot \gamma)\;.
\end{equation}
Also $\scrR^{k)}_\gamma$ is a section of $\Sym^2 TM^*$ along $\gamma$\,.  
However, because of~\eqref{eq:k-ter_Jacobi_on_a_geodesic}, we can define the
$k$-th covariant derivative of the Jacobi operator without reference to a special geodesic as follows:

\bigskip
\begin{definition}
Let $\scrR^{k)}$ denote the section of $\Sym^{k+2}TM^*\otimes \Sym^2 TM^*$ which
is uniquely characterized (via the polarization formula) by
\begin{equation}\label{eq:k-ter_Jacobi}
\scrR^{k)}(\underbrace{\xi,\cdots,\xi}_{k+2};x,y) := \nabla^{k)}_{\xi,\cdots,\xi} R(x,\xi,\xi,y)
\end{equation}
for all $p\in M$ and $\xi,x,y\in T_pM$\,. We will call $\scrR^{k)}$ the
symmetrized $k$-th covariant derivative of the curvature tensor. \end{definition}

For any pseudo-Euclidean vector space $V$ there exist natural inclusions $\Sym^k V^* \supset \Sym^{k-2}V^*\supset \cdots$\,. 
Hence there are also natural inclusions
\[
\Gamma(\Sym^{k+3}TM^*\otimes \Sym^2 TM^*) \supset \Gamma(\Sym^{k+1}TM^*\otimes \Sym^2 TM^*) \supset \cdots
\]
In particular, $\scrR^{k + 1)},\scrR^{k - 1)},\ldots$ may also be seen as sections of $\Sym^{k+3}TM^*\otimes \Sym^2 TM^*$\,, respectively.

\bigskip
\begin{definition}\label{eq:linear_Jacobi_relation}
There exists a linear Jacobi relation of order $k$ if  $\scrR^{k + 1)}$ belongs
to the linear span of $\{\scrR^{k - 1)},\scrR^{k - 3)},\ldots\}$ 
seen as a subspace of the sections of the vector bundle $\Sym^{k + 3}TM^*\otimes\Sym^2TM^*$\,. 
\end{definition} 
In other words, a linear Jacobi relation of order $k$ means that there exist constants $c_{k-1},c_{k-3},\cdots$ such that for all geodesics $\gamma$
\begin{equation}\label{eq:linear_Jacobi_relation_2}
\scrR^{k+1)}_{\gamma} = c_{k-1}\,\g{\dot \gamma}{\dot \gamma} \scrR^{k-1)}_\gamma + c_{k-3}\,\g{\dot \gamma}{\dot \gamma}^2 \scrR^{k-3)}_\gamma + \cdots\;.
\end{equation}


\bigskip
\begin{remark}\label{re:osculating rank}
As was mentioned before, in~\cite{A1,A2,A3,A4,G,GN,NT} a (seemingly) stronger condition called {\em constant Jacobi osculating rank} was considered when the ambient manifold is a Riemannian g.o-space: 
namely, a Riemannian g.o.-space has constant Jacobi osculating rank $k$ if and only if there exists a linear Jacobi relation of order $k$ 
and additionally there exists at least one geodesic $\gamma$ such that  
$\|\dot \gamma(0)\|^2\scrR_\gamma^{k - 1)},\|\dot \gamma(0)\|^4\scrR_\gamma^{k - 3)}(0),\ldots$ are linearly independent. In particular, the linear
Jacobi relation~\eqref{eq:linear_Jacobi_relation_2} has to be minimal then.\footnote{However, an
example of a g.o. space with linear Jacobi relations 
which does not have constant Jacobi osculating does not seem to be known so far.}
\end{remark} 

\section{The Singer invariant}
Let $M$ be a pseudo-Riemannian homogeneous space and $p\in M$\,. The collection $(R,\nabla R,\nabla^{2)}R, \cdots,\nabla^{k)}R)$ will be called the $k$-jet of the curvature tensor. Let $p\in M$ be an arbitrary point and set $V := T_pM$ and let $\so(V)$ denote the  Lie algebra of the orthogonal group. Consider the subalgebra
\begin{equation}\label{eq:k-te_Lie_algebra}
\frakg(k):= \{A\in\so(V) \mid A\cdot R_p = 0\,,\ A\cdot \nabla R_p = 0\,,\ \ldots,\  A\cdot \nabla^{k)} R_p = 0\}
\end{equation}
(here $A\cdot $ means the usual action on arbitrary tensors via algebraic derivation.)

\bigskip
\begin{definition}
Recall from~\cite{GiNi,KT,NTr} that the {\em Singer invariant} of a pseudo-Riemannian homogeneous space is the least number $k_s\in \N_{\geq 0}$  such that
\begin{equation}
\frakg(k_s) = \frakg(k_s + 1)\;.
\end{equation}
\end{definition}
For example, it is known that the Singer invariant of a three- or four-dimensional Riemannian homogeneous space is at most one, cf.~\cite{KT}. For a construction of Riemannian homogeneous spaces with arbitrary high Singer invariant see~\cite{Me} and~\cite[Sec.~6]{Gre}.

Further, let $G$ be a Lie group which acts transitively via isometries on $M$ and $\frakh$ denote the isotropy Lie algebra in $p$\,. Let $\rho_*:\frakh\to\so(V)$ denote the linearized isotropy action. Thus 
\begin{equation}\label{eq:inclusion}
\rho_*(\frakh)\subset \frakg(k_s)\;.
\end{equation}

The Jet Isomorphism Theorem of pseudo-Riemannian geometry (see for
example~\cite[Theorem~3]{J}) immediately implies
that the stabilizer subgroups  in $\SO(V)$ of the $k$-jet of the curvature tensor and
its symmetrized version  $(\scrR_p,\scrR^{1)}_p,\cdots,\scrR^{k)}_p)$ are the same. Hence:

\bigskip
\begin{corollary}
The Lie algebra $\frakg(k)$ has the alternate description
\begin{equation}\label{eq:k-te_Lie_algebra_alternativ}
\frakg(k)= \{A\in\so(V) \mid A\cdot \scrR_p = 0\,,\ A\cdot \scrR^{1)}_p = 0\,,\ \ldots,\  A\cdot \scrR^{k)}_p = 0\}\;,
\end{equation}
the stabilizer Lie algebra of $\jet^k_p \scrR$ in $\so(V)$\,. 
\end{corollary}

Via the description~\eqref{eq:k-te_Lie_algebra_alternativ} of $\frakg(k)$ the following is clear:

\bigskip
\begin{corollary}
Let $M$ be a homogeneous pseudo-Riemannian space with a linear Jacobi
relation of order $k$\,. Then the Singer invariant $k_s$ satisfies
\begin{equation}\label{eq:Abschaetzung}
k_s\leq k\;.
\end{equation}
\end{corollary}

\section{Examples}

\begin{enumerate}
\item
Consider the complex flag-manifold $M^6 := \SU(3)/\rmS(\rmU(1)\times\rmU(1)\times \rmU(1))$ seen as a normal homogeneous space. This is a nearly K{\"a}hler manifold (\cite[pp.142]{BFGK}), in particular Einstein (see also~\cite[p.~577]{WZ}.) According to~\cite{A1}, there is a linear dependence relation
\begin{equation}
\|\dot\gamma\|^4\scrR^{1)}_{\gamma} + 10\|\dot\gamma\|^2\scrR^{3)}_{\gamma}  + 16\,\scrR^{5)}_{\gamma}  = 0
\end{equation}
for every geodesic $\gamma$\,. In fact, the author proves that $M^6$ has
constant Jacobi osculating rang equal to four. 
Hence, it is a priory clear from~\eqref{eq:Abschaetzung} that the Singer invariant is at most four. An explicit calculation with Maple shows that the Singer invariant is actually given by $k_s = 1$\,.
\item Consider the  seven dimensional Berger manifold $V_1= \Sp(2)/\SU(2)$\,, one of the few non-symmetric normal homogeneous spaces with positive sectional curvature. 
This is a strongly irreducible homogeneous space, in particular Einstein.
According to~\cite[Lemma~3.1]{NT} there is a linear dependence relation
\begin{equation}
-\, \|\dot\gamma\|^2\scrR^{1)}_{\gamma}  + \scrR^{3)}_{\gamma}  = 0
\end{equation} 
for every geodesic $\gamma$\,.
In fact, it is shown that $V_1$ has constant Jacobi osculating rank two. Here the Singer invariant is given by $k_s = 0$\,.
\item Consider the  embedding $\rmU(2)\hookrightarrow  \SO(3)\times\SU(3)$ given by the product map $(f,g)$ where $f$ is the adjoint representation on $\su(2)\cong \R^3$  and $g(A):=\diag (A,\det(A)^{-1})$\,. 
The image of $\rmU(2)$ is denoted by $\rmU^\bullet(2)$ and the quotient $V_3 := \SO(3)\times\SU(3)/\rmU^\bullet(2)$ is called the Wilking manifold. 
The Lie algebra $\su(2)\oplus\su(3)$ admits a one-parameter family of positive definite invariant products $\gg_c:= - c\, B_{\su(2)} - B_{\su(3)}$ for $c > 0$ 
inducing different normal homogeneous metrics on $V_3$ all of which have positive sectional curvature (see~\cite{W}). According to~\cite[Lemma~4]{MNT},  for the standard metric (i.e. $c=3/2$)  there is a linear dependence relation
\begin{equation}
-2\, \|\dot\gamma\|^2\scrR^{1)}_{\gamma}  + 5\,\scrR^{3)}_{\gamma}  = 0
\end{equation}
for every geodesic $\gamma$\,. In fact, the authors show that $V_3$ has
constant Jacobi osculating rank two. The Singer invariant is given by
$k_s=0$\,. The Wilking manifold with the standard metric is known to be Einstein (see~\cite{W,Wg}).
\item Consider  {\em Kaplans example} $N^6$\,. This is the simplest H-type group with a two dimensional center. According to a result of~\cite[Lemma~30]{GN} there is a linear dependence relation
\begin{equation}
\|\dot\gamma\|^4\scrR^{1)}_{\gamma}  + 5\, \|\dot\gamma\|^2\scrR^{3)}_{\gamma}  + 4\,\scrR^{5)}_{\gamma}  = 0
\end{equation}
for every geodesic $\gamma$\,. In fact, it is shown that $N^6$ has constant
Jacobi osculating rank four. The Singer invariant satisfies $k_s = 1$\,. 
Note, a non-Abelian nilmanifold is never Einstein (see~\cite[Theorem~2.4]{Mi}.)
\end{enumerate}

\appendix
\section{Some remarks on the Maple calculations of the Singer invariants}
These calculations were done with Maple 17 using the following packages: DifferentialGeometry, LieAlgebra and Tensor.
The structure constants for the Lie algebras
$\fraku(3),\sp(2),\su(3)\oplus\so(3)$ are taken from~\cite{A1,A2,NT,MNT}. 
The calculations provide us with the following additional information on the Lie algebras $\frakg(k)$\,:

\paragraph{For $\SU(3)/\rmS(\rmU(1)\times\rmU(1)\times \rmU(1))$}:
Consider the complex flag-manifold $M^6 := \SU(3)/\rmS(\rmU(1)\times\rmU(1)\times \rmU(1))$ equipped with the normal
homogeneous metric.  According to~\cite{A1}, an orthonormal basis of the six dimensional reductive complement $\frakm$ of the isotropy Lie algebra $\frakh$ is given by
\begin{align}\label{eq:gij_2}
&E_1:= \frac{1}{\sqrt{2}} \left ( \begin{array}{ccc}
0&1&0\\
-1&0&0\\
0&0&0
\end{array}
\right )\;,\ 
&E_2 =  \frac{1}{\sqrt{2}}\left ( \begin{array}{ccc}
0&0&1\\
0&0&0\\
-1&0&0
\end{array}
\right )\;,\ 
&E_3 = \frac{1}{\sqrt{2}}\left ( \begin{array}{ccc}
0&0&0\\
0&0&1\\
0&-1&0
\end{array}
\right )\\
&E_4:=  \frac{1}{\sqrt{2}}\frac{1}{\sqrt{2}} \left ( \begin{array}{ccc}
0&\i&0\\
\i&0&0\\
0&0&0
\end{array}
\right )\;,\ &E_5 =  \frac{1}{\sqrt{2}}\left ( \begin{array}{ccc}
0&0&-\i\\
0&0&0\\
-\i&0&0
\end{array}
\right )\;,\ &E_6 = \frac{1}{\sqrt{2}}\left ( \begin{array}{ccc}
0&0&0\\
0&0&\i\\
0&\i&0
\end{array}
\right )
\end{align}
With respect to this basis, let $\frakt$ denote the usual maximal torus of $\so(6,\R)$ is given by
\[
\{\delta_{1,4}-\delta_{4,1},\delta_{5,2}-\delta_{2,5},\delta_{6,3}-\delta_{3,6}\}_\R
\]
(where $\delta_{ij}$ denotes a matrix with a one in the $i$-th row and the $j$-th column and zeros elsewhere). 
Using Maple, we  calculate that $\frakg(0) = \frakt$ and $\frakg(1) = \rho_*(\frakh)$\,. We conclude from~\eqref{eq:inclusion} that the Singer invariant satisfies $k_s = 1$\,. In fact, already the knowledge of $\frakg(0) = \frakt$
implies $k_s = 1$ since $\SU(3)$ is a simple Lie group which acts effectively on $M^6$ (see~\cite[Theorem~5.2]{WZ}).

\paragraph{For $V_1$ and $V_3$\,:}
Consider the Berger manifold $V_1 := \Sp(2)/\SU(2)$ or the Wilking manifold $V_3 := \SO(3)\times\SU(3)/\rmU^\bullet(2)$ equipped with any normal homogeneous metric. Calculations with Maple show that the Lie algebra $g(0)$ occurring in~\eqref{eq:k-te_Lie_algebra} coincides with $\rho_*(\frakh)$\,. Hence the Singer invariant satisfies $k_s = 0$\,. As a corollary, we see that also for $V_3$ the connected component of the isometry group is equal to $\SO(3)\times\SU(3)$\,.

\paragraph{For  Kaplans example:} Consider the Nilgroup $N^6$ known as  Kaplans example. Following~\cite{D}, we consider the real clifford algebra $\Cl^2_\R:=\Cl(\R^2,\gg)$ associated to the standard inner product $\gg$ on $\R^2$\,. We denote the standard orthonormal basis of $\R^2$ by $\{\i,\j\}$\,. Then $\Cl^2_\R\cong\bbH$ and $\{1,\i,\j,\k:=\i\cdot\j\}$ is a basis of $\bbH$ which satisfies the usual quaternionic relations. 
Further, we consider the standard left-action $\Cl^2_\R\times \bbH\to \bbH, (a,b) \mapsto a\cdot b:= b\,a^{-1}$\,. 
 We set $\frakc:=\{\i,\j\}_\R$ and define a two-step nilpotent Lie bracket on the vector space $\frakn := \bbH\oplus\R^2$ by requiring that $\frakc := \R^2$ 
is the center of $\frakn$ and 
\[
\g{[a,b]}{c} := \g{c\cdot a}{b}
\]
for all $a,b\in\bbH$ and $c\in\frakc$\,.

The corresponding nilpotent Lie group equipped with the so obtained left invariant metric is by definition Kaplans example $N^6$\,. A direct calculation with Maple shows that $k_s\geq 1$ and that the Lie algebra $\frakg(1)$ occuring in~\eqref{eq:k-te_Lie_algebra} is of dimension four (spanned by $\sp(1)$ and the skew-adjoint endomorphism of $\bbH\oplus \R^2$ which is given with respect to the standard basis of $\bbH\oplus \R^2$ by the $6\times 6$ matrix $\delta_{23}-\delta_{32}+\delta_{56}-\delta_{65}$). Moreover, one can show that the latter Lie algebra consists of skew-symmetric derivations of $\frakn$\,. On the other hand, a standard argument shows that the one parameter group generated by a skew-symmetric derivation of $\frakn$ is already a family of isometric Lie group automorphisms of $N$ preserving the identity.
Hence, the isotropy group $H$ is four dimensional and $k_s = 1$\,. This
corresponds to the fact that the isometry group of a nilpotent Lie group equipped with a left invariant metric is ``as small as possible'', see A.~Kaplan~\cite{Ka} and E.~Wilson~\cite{Wi}.

\paragraph{}ACKNOWLEDGMENTS.\ \ \ The author would like to thank Prof.~A.~Naveira for drawing attention to the concept of constant Jacobi osculating rank. I would also like to thank Gregor Weingart for pointing towards the possibility to calculate the Singer invariants with Maple.

\bibliographystyle{amsplain}

\end{document}